\documentclass[11pt,onecolumn,twoside]{IEEEtran} 
\usepackage{amssymb,url,amsmath,amsthm}
\usepackage{here}
\usepackage[pdftex]{graphicx}
\usepackage{bm}
\theoremstyle{theorem}
\newtheorem{theorem}{Theorem}
\newtheorem{proposition}{Proposition}

\newtheorem{lemma}{Lemma}

\newcommand{\eqdef}{:=}
\newcommand{\sq}{\qquad $\square$}

\title{A New Lower Bound for Kullback-Leibler Divergence Based on Hammersley-Chapman-Robbins Bound}
\author{Tomohiro Nishiyama \\ Email: htam0ybboh@gmail.com}

\begin{document} 
\maketitle
\bibliographystyle{plain}
\begin{abstract}
In this paper, we derive a useful lower bound for the Kullback-Leibler divergence (KL-divergence) based on the Hammersley-Chapman-Robbins bound (HCRB).
The HCRB states that the variance of an estimator is bounded from below by the Chi-square divergence and the expectation value of the estimator. By using the relation between the KL-divergence and the Chi-square divergence, we show that the lower bound for the KL-divergence which only depends on the expectation value and the variance of a function we choose. This lower bound can also be derived from an information geometric approach.  Furthermore, we show that the equality holds for the Bernoulli distributions and show that the inequality converges to the Cram\'{e}r-Rao bound when two distributions are very close. We also describe application examples and examples of numerical calculation.
\end{abstract}

\section{Introduction}
The divergences are quantities that measure discrepancy between probability distributions and play key roles in statistics, information theory, signal processing and machine learning. 
For two probability distributions $P$ and $Q$, divergences satisfy the following properties. \\

$D(P\|Q)\geq 0$ with equality if and only if $P=Q$.\\

The Kullback-Leibler divergence (KL-dievergence) \cite{kullback1951information}, the Hellinger distance, and the $\chi^2$-divergence are well-known divergences which are a type of the $\alpha$-divergence \cite{cichocki2010families}.
In particular, the KL-divergence is often used and has good properties that it is invariant under parameter transformations and is compatible with maximum likelihood estimation.

On the other hand, the Hammersley-Chapman-Robbins bound (HCRB) \cite{lehmann2006theory,chapman1951minimum, hammersley1950estimating} states that the variance of an estimator is bounded from below by the $\chi^2 $-divergence and the expected value of the estimator.
From another point of view, the $\chi^2 $-divergence is bounded by the expected value and the variance of a function of interest (FoI). 
Recently, the similar bound for the Hellinger distance was derived in \cite{katsoulakis2017scalable,dashti2016bayesian}

In this paper, we derive a new lower bound for the KL-divergence using the expected value and variance of the FoI.
We show that the equality holds for the Bernoulli distributions and show that we obtain the Cram\'{e}r-Rao bound \cite{cramer1999mathematical, rao1992information} when two distributions are very close. Furthermore, we give an another proof by information geometric approach \cite{amari2016information,amari2010information}. One of the important applications of our inequality is to estimate the lower bound of the KL-divergence between two given data. There are several methods to estimate divergence between given data \cite{nguyen2010estimating, sugiyama2008direct}, our lower bound can be a criterion of validity of estimation results and convergence. Finally, we describe some examples of numerical calculation. \footnote[1]{ \url{https://github.com/nissy220/KL_divergence}}

\section{New bound for the KL-divergence}
\subsection{Main results}
Let us consider two probability distributions (measures) $P$ and $Q$ with the same support $\omega\subseteq\mathbb{R}^n$.
In the following, the integrals is always defined on $\omega$.

The $\alpha$-divergence between $P$ and $Q$ is defined as
\[
  D_\alpha(P\|Q)\eqdef 
\begin{cases}
    \frac{1}{\alpha(\alpha-1)}\biggl(\int p(x)^\alpha q(x)^{1-\alpha} d\lambda(x)-1\biggr) & (\alpha\neq 0,1) \\
    \int q(x)\log \frac{q(x)}{p(x)} d\lambda(x) & (\alpha=0) \\
    \int p(x)\log \frac{p(x)}{q(x)} d\lambda(x) & (\alpha=1),
  \end{cases}
\]
where $p$ and $q$ denote the Radon-Nikodym derivatives of probability measures $P$ and $Q$ respectively and $\lambda$ denotes the Lebesgue measure or the counting measure.

The KL-divergence, the squared Hellinger distance and the $\chi^2$-divergence are the cases for $\alpha=1, \frac{1}{2},2$.
\begin{align}
\label{def_kl}
\mathrm{KL}(P\|Q)\eqdef \int p(x)\log\frac{p(x)}{q(x)} d\lambda(x) =D_1(P\|Q)\\
\mathrm{Hel}^2(P,Q)\eqdef \int (\sqrt{q(x)}-\sqrt{p(x)})^2 d\lambda(x)=\frac{1}{2}D_\frac{1}{2}(P\|Q)\\
\chi^2(P\|Q)\eqdef \int \frac{(q(x)-p(x))^2}{q(x)} d\lambda(x) = 2D_{2}(P\|Q). 
\end{align}

Next, we consider a function of interest (FoI) $f: \omega\rightarrow\mathbb{R}$ and consider an expected value and a variance of the FoI for $P$.
These quantities are defined as $E_P[f]\eqdef \int p(x)f(x) d\lambda(x)$ and $\mathrm{Var}_P(f)\eqdef E_P[f^2]-E_P[f]^2$. 

The $\chi^2$-divergence and the squared Hellinger distance satisfy the following inequalities.
\begin{align}
\label{ineq_hcrb}
\chi^2(P\|Q) \geq \frac{(E_Q-E_P)^2}{V_Q} \\ 
\mathrm{Hel}^2(P,Q) \geq \frac{(E_Q-E_P)^2}{2\bigl(V_P+V_Q+\frac{1}{2}(E_Q-E_P)^2\bigr)},
\end{align}
where $E_P, E_Q$ and $V_P, V_Q$ denote $E_P[f], E_Q[f]$ and $\mathrm{Var}_P(f), \mathrm{Var}_Q(f)$, respectively.
The inequality (\ref{ineq_hcrb}) is the HCRB for general probability densities.

The main purpose of this section is to derive a similar lower bound which only depends on $E_P, E_Q, V_P, V_Q$ for the KL-divergence.
We first show the result.
\begin{theorem} [New lower bound for the KL-divergence]
\label{th_LB}
Let $P$ and $Q$ are probability distributions with the same support $\omega$.

Let $f:\omega \rightarrow \mathbb{R}$ be a FoI with finite $E_P, E_Q, V_P, V_Q$ and $V_P, V_Q > 0$, where $E_P, E_Q$ and $V_P, V_Q$ denote $E_P[f], E_Q[f]$ and $\mathrm{Var}_P(f), \mathrm{Var}_Q(f)$, respectively.

Then,
\begin{align}
\label{th_eq1}
\mathrm{KL}(P\|Q) \geq \frac{(A-2V_P)}{D} \tanh^{-1} \frac{D}{A}+\frac{1}{2}\log\frac{V_P}{V_Q},
\end{align}
where  $A\eqdef (E_Q-E_P)^2 + V_P+V_Q$, $D\eqdef \sqrt{A^2-4V_PV_Q}$.
The equality holds if and only if there exists a function $C(t)$ which satisfies 
\begin{align}
\frac{q(x)-p(x)}{p(x)+t(q(x)-p(x))}=\bigl(f(x)-E_{P}-t(E_{Q}-E_{P})\bigr)C(t)
\end{align}
for all $x\in\omega$ and $t\in[0,1]$, where $p$ and $q$ denote the Radon-Nikodym derivatives of $P$ and $Q$, respectively.
\end{theorem}

The basic idea to show this inequality is to use a relation between $D_\alpha$ and $D_{\alpha+1}$. Considering a mixture distribution $r(x;t)=p(x)+t(q(x)-p(x))$ for $t\in[0,1]$, the derivative of  $D_\alpha(P\|R(t))$ with respect to $t$ can be represented using $D_\alpha(P\|Q)$ and $D_{\alpha+1}(P\|Q)$ as shown later.
Since the KL-divergence and the $\chi^2$-divergence are the case for $\alpha=1$ and $\alpha=2$ respectively, by applying the HCRB and integrating from 0 to 1 with respect to $t$, we can derive the inequality. We show details below.

\begin{lemma}
\label{lem_alpha}
Let $P$ and $Q$ are probability distributions with the same support $\omega$.

Let $r(x;t)$ be the Radon-Nikodym derivative of a probability distribution $R(t)$ and $r(x;t)\eqdef p(x)+t(q(x)-p(x))$ for $t\in[0,1]$.

Then,
\begin{align}
\frac{d}{dt}D_\alpha(P\|R(t))=\frac{(1-\alpha)}{t}D_\alpha(P\|R(t))+\frac{(1+\alpha)}{t}D_{\alpha+1}(P\|R(t)).
\end{align}
\end{lemma}
\noindent\textbf{Proof}.
For $\alpha\neq 0,1$,
\begin{align}
\label{lem1_eq1}
\frac{d}{dt}D_\alpha(P\|R(t))&=-\frac{1}{\alpha}\int (q(x)-p(x))p(x)^\alpha r(x;t)^{-\alpha} d\lambda(x)
=-\frac{1}{\alpha}\int \frac{(r(x;t)-p(x))}{t} p(x)^\alpha r(x;t)^{-\alpha} d\lambda(x) \\
&=\frac{(1-\alpha)}{t}D_\alpha(P\|R(t))+\frac{(1+\alpha)}{t}D_{\alpha+1}(P\|R(t)).
\end{align}

For $\alpha=1$, 
\begin{align}
\frac{d}{dt}D_1(P\|R(t))&=-\int\frac{p(x)(q(x)-p(x))}{r(x;t)}d\lambda(x) \\
&=\frac{1}{t}\int\frac{p(x)(p(x)-r(x;t))}{r(x;t)}d\lambda(x)=\frac{1}{t}\int\frac{(r(x;t)-p(x))^2}{r(x;t)}d\lambda=\frac{2}{t}D_2(P\|R(t)).
\end{align}
We can show the inequality for $\alpha=0$ in the same way.
Hence, the result follows. \sq
\begin{lemma}  [HCRB]
\label{lem_chapman}
Let $P$ and $Q$ are probability distributions with the same support $\omega$.

Let $f:\omega \rightarrow \mathbb{R}$ be a FoI with finite $E_P, E_Q, V_P, V_Q$ and $V_Q > 0$. 

Then, 
\begin{align}
\label{lem2_eq1}
\chi^2(P\|Q) \geq \frac{(E_Q-E_P)^2}{V_Q},
\end{align}
with equality if and only if there exists a constant $C$ which satisfies  $\frac{q(x)-p(x)}{q(x)}=C(f(x)-E_Q)$ for all $x\in\omega$.
\end{lemma}
\noindent\textbf{Proof.}
Consider the following quantity.
\begin{align}
E_Q-E_P=\int (f(x)-E_Q) (q(x)-p(x)) d\lambda(x)=\int \sqrt{q(x)} (f(x)-E_Q) \frac{(q(x)-p(x))}{\sqrt{q(x)}} d\lambda(x).
\end{align}
Applying the Cauchy-Schwartz inequality yields
\begin{align}
(E_Q-E_P)^2\leq V_Q\int \frac{(q(x)-p(x))^2}{q(x)} d\lambda(x) =V_Q\chi^2(P\|Q) .
\end{align}
Since $V_Q>0$, we obtain the desired inequality.
Since the FoI isn't a constant from the assumption  $V_Q>0$, the equality holds if and only if there exists a constant $C$ which satisfies $\frac{q(x)-p(x)}{\sqrt{q(x)}}=C\sqrt{q(x)} (f(x)-E_Q)$ for all $x\in\omega$.
Hence, the result follows. \sq

\noindent\textbf{Proof of Theorem 1}.
Applying Lemma \ref{lem_alpha} for $\alpha=1$ yields
\begin{align}
\label{th_eq1}
\frac{d}{dt}\mathrm{KL}(P\|R(t))= \frac{\chi^2(P\|R(t))}{t}.
\end{align}
Applying Lemma \ref{lem_chapman}, we have
\begin{align}
\label{th_eq2}
\frac{d}{dt}\mathrm{KL}(P\|R(t))\geq \frac{(E_{R(t)}-E_P)^2}{t\mathrm{V}_{R(t)}}.
\end{align}
Since $E_{R(t)}[f]=tE_Q[f]+(1-t)E_P[f]$ and $E_{R(t)}(f^2)=tE_Q[f^2]+(1-t)E_P[f^2]$, we have 
$\mathrm{Var}_{R(t)}(f)=E_{R(t)}(f^2)-E_{R(t)}[f]^2=tV_Q+(1-t)V_P+t(1-t)(E_Q-E_P)^2$.
Substituting this equality into (\ref{th_eq2}), we have
\begin{align}
\label{th_eq3}
\frac{d}{dt}\mathrm{KL}(P\|R(t))\geq \frac{t(E_Q-E_P)^2}{t(1-t)(E_Q-E_P)^2+(1-t)V_P+tV_Q}.
\end{align}
Integrating from 0 to 1 with respect to $t$, we have
\begin{align}
\label{th_eq4}
\mathrm{KL}(P\|Q)\geq \int_{0}^1 \frac{t(E_Q-E_P)^2}{t(1-t)(E_Q-E_P)^2+(1-t)V_P+tV_Q}dt.
\end{align}

Using the following formula
\begin{align}
\label{th_eq5}
\int_{0}^x \frac{at}{-at^2+(a+b)t+c} dt \\ \nonumber
=\frac{a+b}{\sqrt{a^2+2a(b+2c)+b^2}}\tanh^{-1}\biggl(\frac{a(2x-1)-b}{\sqrt{a^2+2a(b+2c)+b^2}}\biggr)-\frac{1}{2}\log(x(-ax+a+b)+c)+\mathrm{const.}
\end{align}
and putting $a=(E_Q-E_P)^2$, $b=V_Q-V_P$, $c=V_P$ and using the definition $A=(E_Q-E_P)^2+V_P+V_Q$ and $D=\sqrt{A^2-4V_PV_Q}$, we have
\begin{align} 
\label{th_eq6}
\int_{0}^1 \frac{t(E_Q-E_P)^2}{t(1-t)(E_Q-E_P)^2+(1-t)V_P+tV_Q} dt \\ \nonumber
=\frac{(A-2V_P)}{D}\biggl(\tanh^{-1}\biggl(\frac{a-b}{D}\biggr) + \tanh^{-1}\biggl(\frac{a+b}{D}\biggr)\biggr)-\frac{1}{2}\log\frac{V_Q}{V_P}.
\end{align}
By applying the addition formula of the inverse hyperbolic function $\tanh^{-1}u +\tanh^{-1}v=\tanh^{-1}\bigl(\frac{u+v}{1+uv}\bigr) $ and $D^2+a^2-b^2=2aA$, we have 
\begin{align}
\label{th_eq7}
\int_{0}^1 \frac{t(E_Q-E_P)^2}{t(1-t)(E_Q-E_P)^2+(1-t)V_P+tV_Q} dt=
\frac{(A-2V_P)}{D}\tanh^{-1}\biggl(\frac{D}{A}\biggr)+\frac{1}{2}\log\frac{V_P}{V_Q}.
\end{align}
By combining with (\ref{th_eq4}), we obtain the desired inequality.
From Lemma \ref{lem_chapman}, the equality holds if there exists a function $\tilde{C}(t)$ which satisfies
\begin{align}
\label{th_eq7}
\frac{r(x;t)-p(x)}{r(t;x)}=\tilde{C}(t)(f(x)-E_{R(t)})
\end{align}
for all $x\in\omega$ and $t\in[0,1]$.   
By substituting $r(x;t)=p(x)+t(q(x)-p(x))$ and $E_{R(t)}=E_P+t(E_Q-E_P)$ into this equality and putting $C(t)=\frac{\tilde{C}(t)}{t}$, we obtain the equality condition.   \sq

\begin{proposition}
When $\omega=\{0,1\}$ and the FoI is $x$, the Bernoulli distributions satisfy the equality condition of Theorem \ref{th_LB}.
\end{proposition}
\noindent\textbf{Proof}.
Let $P(x)$ and $Q(x)$ be probability mass functions of the Bernoulli distributions and $P(0)=1-p, P(1)=p$ and $Q(0)=1-q, Q(1)=q$.
We choose a function $C(t)$ which satisfies
\begin{align}
\frac{Q(0)-P(0)}{P(0)+t(Q(0)-P(0))}=-(E_P+t(E_Q-E_P))C(t).
\end{align}
This is the equality condition of Theorem \ref{th_LB} for $x=0$.
Since $E_P+t(E_Q-E_P)=P(1)+t(Q(1)-P(1))$, $P(0)+t(Q(0)-P(0))=1-(E_P+t(E_Q-E_P))$ and $Q(0)-P(0)=-(Q(1)-P(1))$, we have
\begin{align}
\frac{Q(1)-P(1)}{P(1)+t(Q(1)-P(1))}=(1-(E_P+t(E_Q-E_P))C(t).
\end{align}
This is the equality condition for $x=1$.
Hence, the result follows. \sq

\begin{proposition}
For a parameter $\theta\in\mathbb{R}$, let $p(x;\theta)$ and $p(x;\theta+\delta\theta)$ be probability density functions of probability distributions $P$ and $Q$, respectively.
When $\delta\theta\rightarrow 0$, the inequality of Theorem \ref{th_LB} converges to the Cram\'{e}r-Rao bound.

The Cram\'{e}r-Rao bound is 
\begin{align}
\label{CR_bound}
\mathrm{Var}_P(f)\geq \frac{\psi'(\theta)^2}{I(p)},
\end{align}
where $I(p)\eqdef \int \frac{p'(x,\theta)^2}{p(x,\theta)} d\lambda(x)$ is the Fisher information, $\psi(\theta)\eqdef E_P[f](\theta)$ and $'$ denotes the derivative with respect to $\theta$.
\end{proposition}
\noindent\textbf{Proof}.
Substituting $E_Q-E_P=(E_P[f])'(\theta)\delta\theta+O(\delta\theta^2)$ and $V_P=V_Q+O(\delta\theta)$ into (\ref{th_eq4}) yields
\begin{align}
\mathrm{KL}(P\|Q)\geq \delta\theta^2\int_{0}^1 \frac{t\psi'(\theta)^2}{V_P}dt+O(\delta\theta^3)= \delta\theta^2\frac{\psi'(\theta)^2}{2\mathrm{Var}_P(f)}+O(\delta\theta^3).
\end{align}
On the other hand, since $\mathrm{KL}(P\|Q)=\frac{1}{2}I(p)\delta\theta^2+O(\delta\theta^3)$ holds (see \cite{kullback1951information}), in the limit $\delta\theta\rightarrow 0$, we have 
\begin{align}
I(p)\geq \frac{\psi'(\theta)^2}{\mathrm{Var}_P(f)}.
\end{align}
Hence, the result follows. \sq

One of the important application of our lower bound is divergence estimation between two given data.
Since we can approximately calculate the expected value of the FoI of the data $\{x_i\}$ such as 
$E[f]=\frac{1}{N}\sum_i f(x_i) $ and $\mathrm{Var}(f)=\frac{1}{N}\sum_i f(x_i)^2 -E[f]^2 $, we can easily calculate the lower bound of the KL-divergence by using Theorem \ref{th_LB} even if the true distributions is unknown. Hence, Our lower bound can be a criteria to judge the validity of estimation results and to abort the iteration.

\subsection{Another proof of Theorem 1}
In this subsection, we prove Theorem \ref{th_LB} by the information geometric approach.
First, we derive the relation between the derivative of the canonical divergence along the geodesic on the dually flat manifold and the Fisher information metric.
Next, by applying the Cram\'{e}r-Rao bound to this relation, we prove Theorem \ref{th_LB}.

Let $P_i (i=0,1,\dots ,d)$ be probability distributions (measures) on the same support $\omega\subseteq\mathbb{R}^n$ and $S\eqdef\{\boldsymbol\eta|\boldsymbol\eta=(\eta_1, \eta_2, \dots, \eta_d)\in \mathbb{R}^d, \eta_i\geq 0 ,\sum_{i=1}^d \eta_i \leq 1\}$.
Consider mixture distributions $r(x;\boldsymbol\eta)\eqdef (1-\sum_{i=1}^d \eta_i)p_0(x)+\sum_{i=1}^d \eta_i p_i(x) $ for $\boldsymbol\eta\in S$ and the function called potential $F(\boldsymbol\eta)\eqdef E_{r(x;\boldsymbol\eta)}[\log r(x;\boldsymbol\eta)]$, where $p_i$ denote the Radon-Nikodym derivatives of probability measures $P_i$ and $E_{r(x;\boldsymbol\eta)}$ denotes the expected value of $r(x;\boldsymbol\eta)$.
By the definition of $F$, $F$ is the convex function and the negative Shannon entropy.
The potential $F$ and the coordinate $\boldsymbol\eta$ define the dually flat structure in the manifold $M=\{r(x;\boldsymbol\eta)|\boldsymbol\eta\in S \}$. 
We can introduce the canonical Bregman divergence \cite{nielsen2018elementary} on $M$.
\begin{align}
D(\boldsymbol\eta_1\|\boldsymbol\eta_2)\eqdef F(\boldsymbol\eta_1)-F(\boldsymbol\eta_2)-\sum_i \partial^iF(\boldsymbol\eta_2)(\eta_{1i}-\eta_{2i}),
\end{align}
where $\partial^i\eqdef \frac{\partial}{\partial\eta_i}$. 
The Riemannian metric $g^{ij}$ on $M$ is $g^{ij}=\partial^i\partial^j F$ and $g^{ij}$ is equal to the Fisher information metric $I^{ij}$.
\begin{align}
\label{eq_Riemann_Fisher}
g^{ij}(\boldsymbol\eta)=\partial^i\partial^j F(\boldsymbol\eta)=I^{ij}(\boldsymbol\eta)\eqdef E_{r(x;\boldsymbol\eta)}[\partial^i l_{\boldsymbol\eta}\partial^j l_{\boldsymbol\eta}],
\end{align}
where $l_{\boldsymbol\eta}$ denotes $\log r(x;\boldsymbol\eta)$.
We can easily confirm this equality by the definition of the potential $F$ and $r(x;\boldsymbol\eta)$.
In the same way, we have
\begin{align}
\label{eq_canonical_KL}
D(\boldsymbol\eta_1\|\boldsymbol\eta_2)=\mathrm{KL}(R(\boldsymbol\eta_1)\|R(\boldsymbol\eta_2)),
\end{align}
where $R(\boldsymbol\eta)$ denotes probability distributions which correspond to $r(x;\boldsymbol\eta)$.
The geodesic on $M$ can be written as $\boldsymbol\eta(t)=\bm{a}t+\bm{b}$ for the constant vector $\bm{a}, \bm{b}\in \mathbb{R}^d$ and a parameter $t\in \mathbb{R}$.

From the definition of the canonical divergence, the derivative of the canonical divergence along the geodesic is 
\begin{align}
\label{eq_derivative_Fisher_multi}
\frac{d}{dt}D(\boldsymbol\eta(0)\|\boldsymbol\eta(t))=t\sum_{i,j}\partial^i\partial^j F(\boldsymbol\eta(t))a_ia_j=t\sum_{i,j}g^{ij}(\boldsymbol\eta(t))a_ia_j=t\sum_{i,j}I^{ij}(\boldsymbol\eta(t))a_ia_j,
\end{align}
where we use (\ref{eq_Riemann_Fisher}).
This equality represents the relation between the derivative of the canonical divergence along the geodesic and the Fisher information metric.

In the case $M$ is 1-dimensional and $\eta_1(t)=t$, (\ref{eq_derivative_Fisher_multi}) is simplified as 
\begin{align}
\label{eq_derivative_Fisher}
\frac{d}{dt}D(0\|t)=tI(t).
\end{align}
By putting $p_0(x)=p(x)$, $p_1(x)=q(x)$, $r(x;t)=p(x)+t(q(x)-p(x))$ and applying the Cram\'{e}r-Rao bound (\ref{CR_bound}), we have
\begin{align}
\label{CR_bound_2}
V_{R(t)}=\mathrm{Var}_{R(t)}(f)\geq \frac{(E_Q-E_P)^2}{I(t)},
\end{align}
where we use $\psi(t)=E_{R(t)}[f]=(E_Q-E_P)t+E_P$.
From (\ref{eq_canonical_KL}), (\ref{eq_derivative_Fisher}) and (\ref{CR_bound_2}), we have
\begin{align}
\frac{d}{dt}\mathrm{KL}(P\|R(t))\geq \frac{t(E_Q-E_P)^2}{V_{R(t)}}.
\end{align}
By using $\mathrm{Var}_{R(t)}(f)=E_{R(t)}(f^2)-E_{R(t)}[f]^2=tV_Q+(1-t)V_P+t(1-t)(E_Q-E_P)^2$, we have the same inequality as (\ref{th_eq3}).
The rest of the proof is the same as the previous subsection.

\section{Examples}
We show four examples using Theorem \ref{th_LB}.
\subsection{Simple discrete distribution}
We consider the case $\omega=\{1,2,3,4\}$, $P(1)=P(2)=P(3)=P(4)=0.25$ and $Q(1) = 0.1, Q(2)=0.2, Q(3)=0.3, Q(4) = 0.4$.
Then, we have $E_P=2.5$, $E_Q=3.0$ and $V_P=1.25, V_Q=1.0$. 
The KL-divergence is equal to $0.121777274287$ and our lower bound is equal to $0.111571775657$.
\subsection{Normal distribution}
The probability density function of the normal distribution is 
\begin{align}
p(x;\mu, \sigma)=\frac{1}{\sqrt{2\pi}\sigma}\exp\biggl(-\frac{(x-\mu)^2}{2\sigma^2}\biggr). 
\end{align}
and the KL-divergence is 
\begin{align}
KL(\mu_P, \sigma_P\| \mu_Q, \sigma_Q)=\frac{(\mu_Q-\mu_P)^2}{2\sigma_Q^2}+\frac{\sigma_P^2}{2\sigma_Q^2}-\frac{1}{2}+\log\frac{\sigma_Q}{\sigma_P}.
\end{align}

\begin{figure}[H]
 \begin{center}
  \includegraphics[width=100mm, height = 50mm]{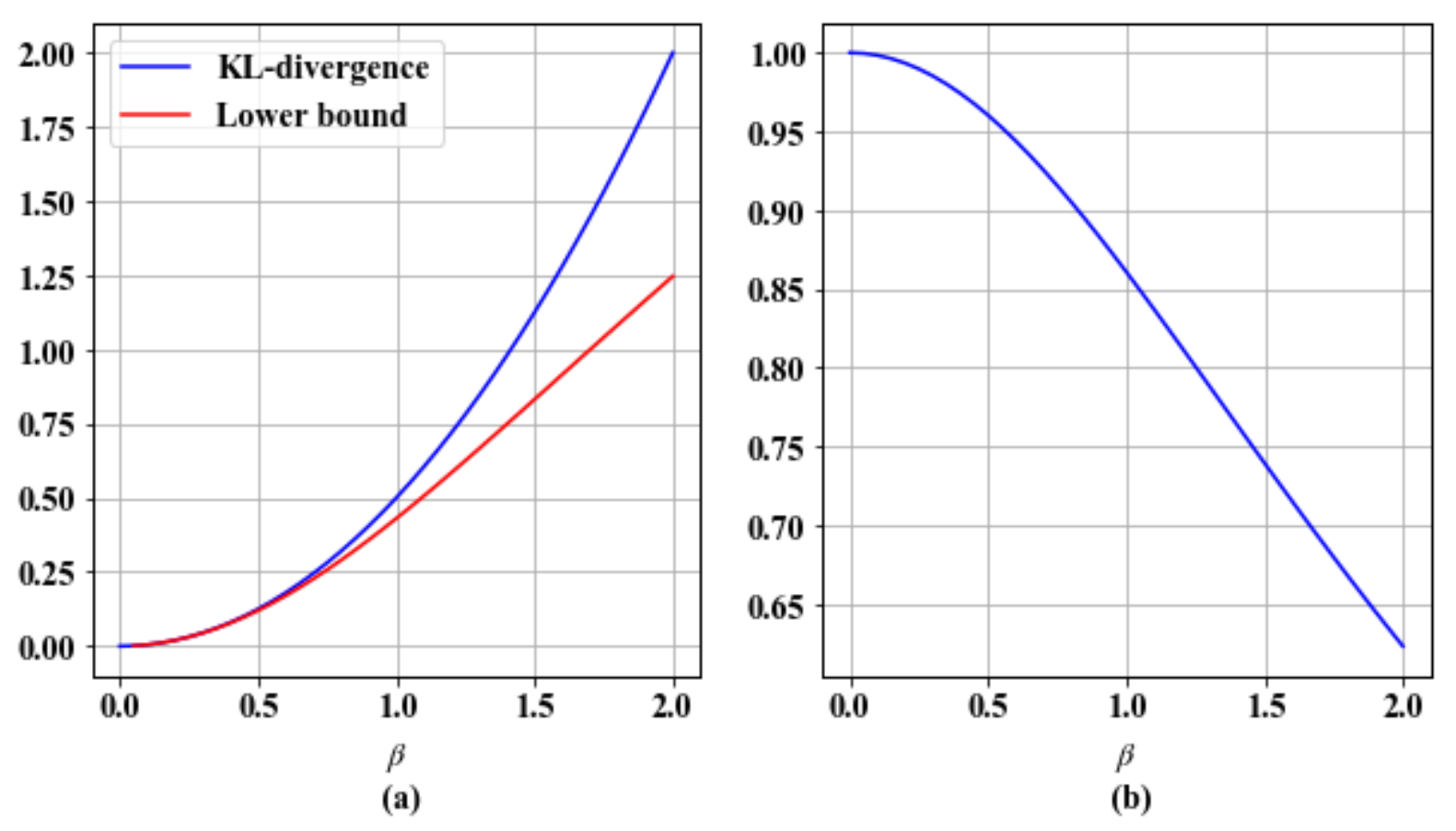}
 \end{center}
 \caption{(a) Comparison of the KL-divergence (blue line) and our lower bound (red line) for the normal distribution. (b) The ratio of our lower bound to the KL-divergence.  }
 \label{fig:one}
\end{figure}

1) We consider the case $\sigma_P=\sigma_Q=\sigma$.
In this case, the KL-divergence is simplified as $KL(\mu_P, \sigma_P\| \mu_Q, \sigma_Q)=\frac{(\mu_Q-\mu_P)^2}{2\sigma^2}$.
Choosing the foI as $f(x)=x$, we have $E_P=\mu_P$, $E_Q=\mu_Q$ and $V_P=V_Q=\sigma^2$. 
Figure \ref{fig:one}. displays the result of comparison of the KL-divergence and our lower bound for $\beta\eqdef \frac{\mu_Q-\mu_P}{\sigma}$.

\begin{figure}[H]
 \begin{center}
  \includegraphics[width=100mm, height = 50mm]{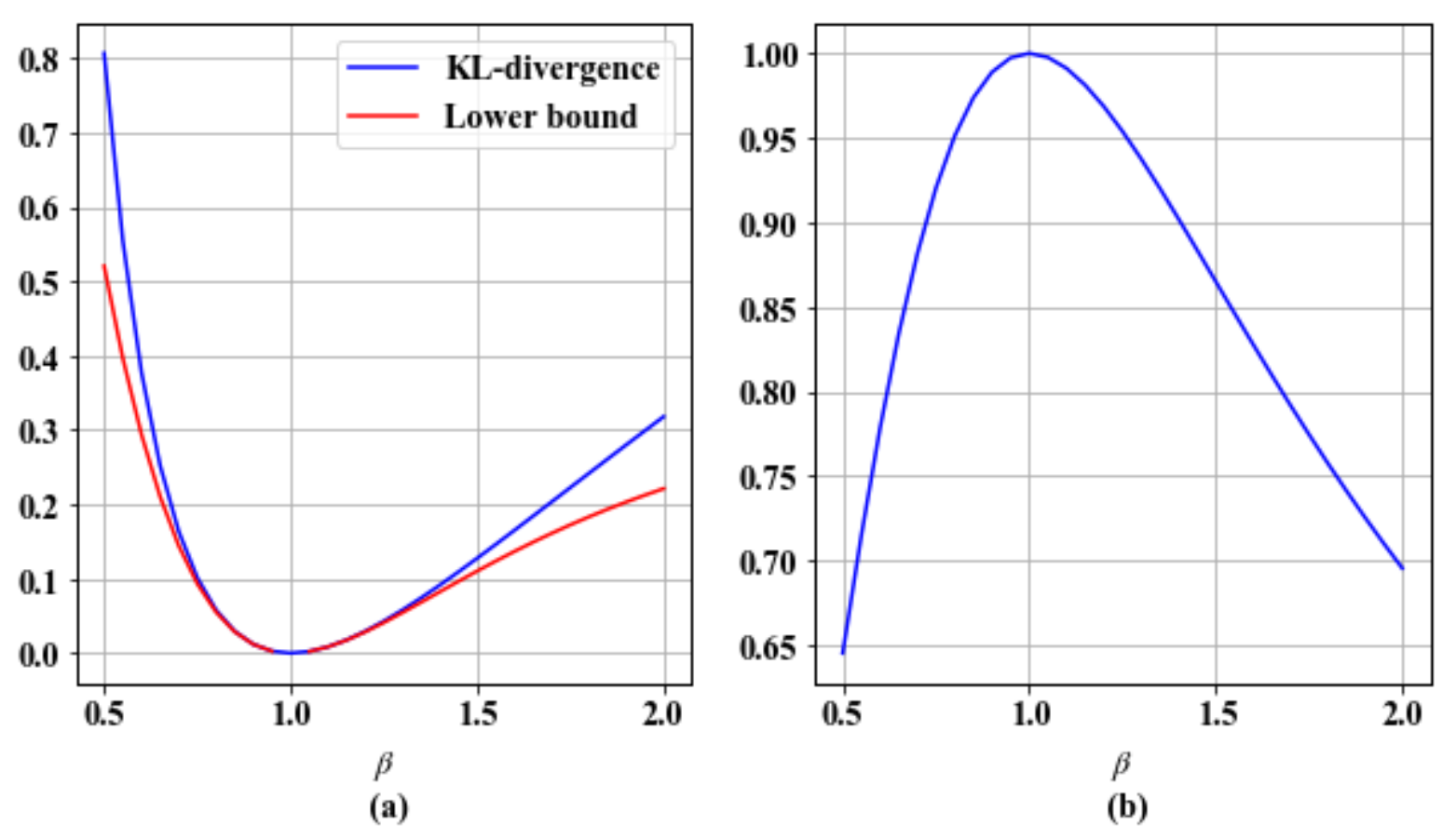}
 \end{center}
 \caption{(a) Comparison of the KL-divergence (blue line) and our lower bound (red line) for the normal distribution. (b) The ratio of our lower bound to the KL-divergence. }
 \label{fig:two}
\end{figure}

2) We consider the case $\mu_P=\mu_Q=0$.
In this case, the KL-divergence is simplified as $\frac{\sigma_P^2}{2\sigma_Q^2}-\frac{1}{2}+\log\frac{\sigma_Q}{\sigma_P}$.
If we choose the foI as $f(x)=x$, the lower bound is trivial from (\ref{th_eq4}).
Hence, we choose the foI as $f(x)=x^2$ and we have $E_P=\sigma_P^2$, $E_Q=\sigma_Q^2$ and $V_P=2\sigma_P^4$ $V_Q=2\sigma_Q^4$.
Figure \ref{fig:two}. displays the result of comparison of the KL-divergence and our lower bound for $\beta\eqdef \frac{\sigma_Q}{\sigma_P}$.

\subsection{Exponential distribution}
\begin{figure}[H]
 \begin{center}
  \includegraphics[width=100mm, height = 50mm]{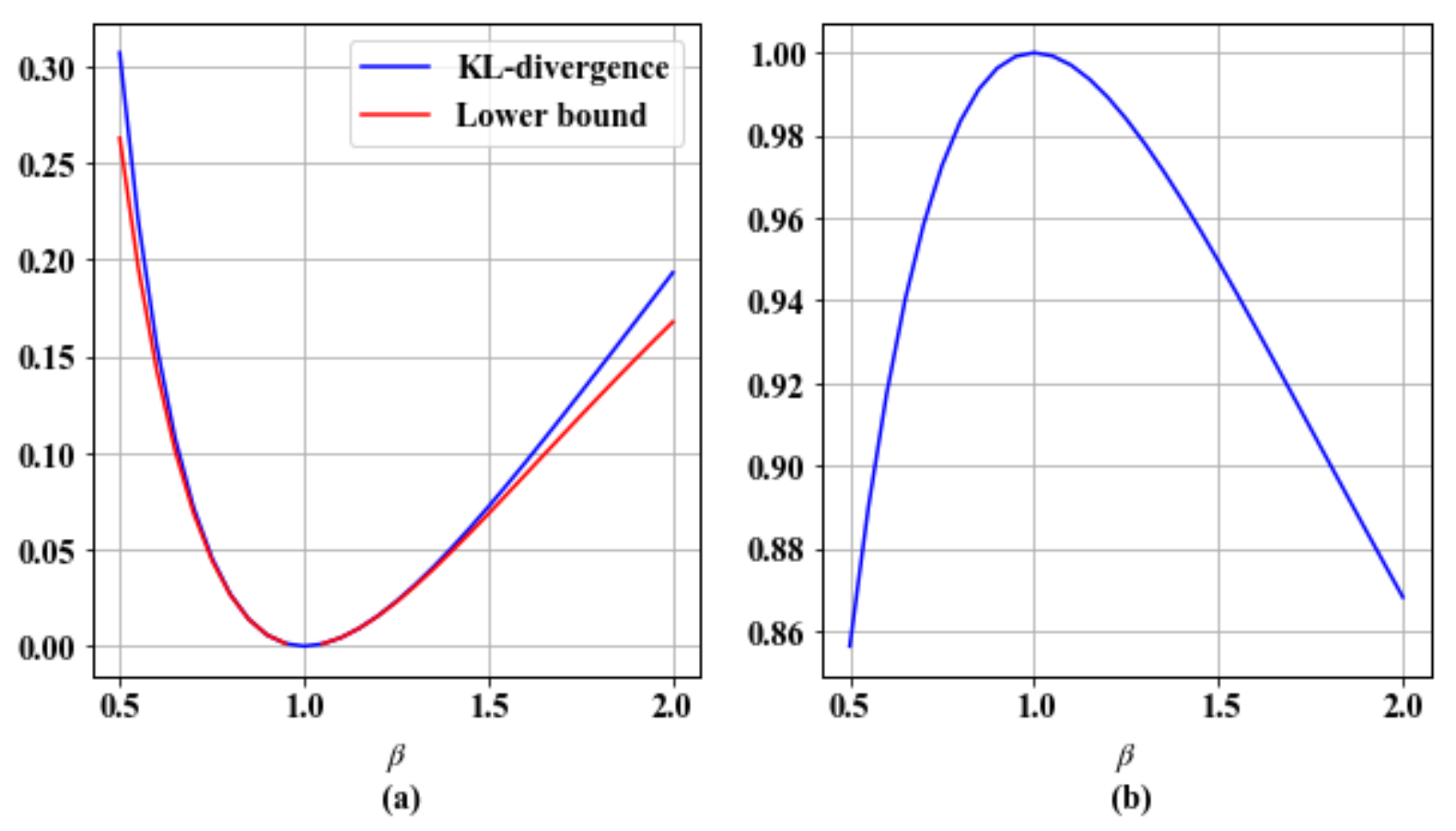}
 \end{center}
 \caption{(a) Comparison of the KL-divergence (blue line) and our lower bound (red line) for the exponential distribution. (b) The ratio of our lower bound to the KL-divergence. }
 \label{fig:three}
\end{figure}

The probability density function of the exponential distribution is 
\begin{align}
p(x;\nu)=\frac{1}{\nu}\exp\bigl(-\frac{1}{\nu} x\bigr).
\end{align}
and the KL-divergence is 
\begin{align}
KL(\nu_P\| \nu_Q)=\frac{\nu_P}{\nu_Q}-1+\log\frac{\nu_Q}{\nu_P}.
\end{align}
Choosing the foI as $f(x)=x$, we have $E_P=\nu_P$, $E_Q=\nu_Q$ and $V_P=\nu_P^2$,  $V_Q=\nu_Q^2$. 
Figure \ref{fig:three}. displays the result of comparison of the KL-divergence and our lower bound for $\beta\eqdef \frac{\nu_Q}{\nu_P}$.

From these examples, we can confirm that our lower bounds behave similarly to the actual value of the KL-divergence.
 
\section{Conclusion}
We have derived the new lower bound for the KL-divergence which only depends on the expected value and the variance of the FoI by applying the HCRB and we have given another proof by information geometric approach.
We also have shown that the equality holds for the Bernoulli distributions and we have obtained the Cram\'{e}r-Rao bound from our lower bound when two distributions are very close.
Furthermore, we have described that our bound behaves similarly to the KL-divergence by some examples of numerical calculation.
One of the important point of our bound is easiness of calculation even if the true distribution is unknown.
We expect the application range of it will expand.
\bibliography{reference_KL}

\begin{thebibliography}{10}

\bibitem{amari2016information}
Shun-ichi Amari.
\newblock {\em Information geometry and its applications}.
\newblock Springer, 2016.

\bibitem{amari2010information}
Shun-ichi Amari and Andrzej Cichocki.
\newblock Information geometry of divergence functions.
\newblock {\em Bulletin of the Polish Academy of Sciences: Technical Sciences},
  58(1):183--195, 2010.

\bibitem{chapman1951minimum}
Douglas~G Chapman, Herbert Robbins, et~al.
\newblock Minimum variance estimation without regularity assumptions.
\newblock {\em The Annals of Mathematical Statistics}, 22(4):581--586, 1951.

\bibitem{cichocki2010families}
Andrzej Cichocki and Shun-ichi Amari.
\newblock Families of alpha-beta-and gamma-divergences: Flexible and robust
  measures of similarities.
\newblock {\em Entropy}, 12(6):1532--1568, 2010.

\bibitem{cramer1999mathematical}
Harald Cram{\'e}r.
\newblock {\em Mathematical methods of statistics}, volume~9.
\newblock Princeton university press, 1999.

\bibitem{dashti2016bayesian}
Masoumeh Dashti and Andrew~M Stuart.
\newblock The bayesian approach to inverse problems.
\newblock {\em Handbook of Uncertainty Quantification}, pages 1--118, 2016.

\bibitem{hammersley1950estimating}
JM~Hammersley.
\newblock On estimating restricted parameters, jr statist.
\newblock {\em Soc.(B)}, 12, 1950.

\bibitem{katsoulakis2017scalable}
Markos~A Katsoulakis, Luc Rey-Bellet, and Jie Wang.
\newblock Scalable information inequalities for uncertainty quantification.
\newblock {\em Journal of Computational Physics}, 336:513--545, 2017.

\bibitem{kullback1951information}
Solomon Kullback and Richard~A Leibler.
\newblock On information and sufficiency.
\newblock {\em The annals of mathematical statistics}, 22(1):79--86, 1951.

\bibitem{lehmann2006theory}
Erich~L Lehmann and George Casella.
\newblock {\em Theory of point estimation}.
\newblock Springer Science \& Business Media, 2006.

\bibitem{nguyen2010estimating}
XuanLong Nguyen, Martin~J Wainwright, and Michael~I Jordan.
\newblock Estimating divergence functionals and the likelihood ratio by convex
  risk minimization.
\newblock {\em IEEE Transactions on Information Theory}, 56(11):5847--5861,
  2010.

\bibitem{nielsen2018elementary}
Frank Nielsen.
\newblock An elementary introduction to information geometry.
\newblock {\em arXiv preprint arXiv:1808.08271}, 2018.

\bibitem{rao1992information}
C~Radhakrishna Rao.
\newblock Information and the accuracy attainable in the estimation of
  statistical parameters.
\newblock In {\em Breakthroughs in statistics}, pages 235--247. Springer, 1992.

\bibitem{sugiyama2008direct}
Masashi Sugiyama, Taiji Suzuki, Shinichi Nakajima, Hisashi Kashima, Paul von
  B{\"u}nau, and Motoaki Kawanabe.
\newblock Direct importance estimation for covariate shift adaptation.
\newblock {\em Annals of the Institute of Statistical Mathematics},
  60(4):699--746, 2008.

\end{thebibliography}
\end{document}